\documentclass{article}      
\usepackage{amssymb}
\usepackage{amsmath}
\usepackage{latexsym}
\usepackage[mathscr]{euscript}
\usepackage{graphicx}
\usepackage{amscd}
\usepackage{amsthm} 
\usepackage{fullpage} 
\usepackage{wrapfig}
\newtheorem{theorem}{Theorem}

\newtheorem{definition}{Definition}

\newcommand{\bz}{\mathbb{Z}}
\newcommand{\br}{\mathbb{R}}

\newcommand{\p}{\partial}
\newcommand{\cc}{\mathcal{C}}
\newcommand{\ce}{\mathcal{E}}
\newcommand{\cj}{\mathcal{J}}

\newcommand{\cl}{\mathcal{L}}

\newcommand{\cb}{\mathcal{B}}

\newcommand{\ch}{\mathcal{H}}
\newcommand{\cm}{\mathcal{M}}

\newcommand{\cy}{\mathcal{Y}}
\newcommand{\cdd}{\mathcal{D}}

\newcommand{\ca}{\mathcal{A}}

\newcommand{\hk}{\hookrightarrow}

\newcommand{\med}{\medskip}
\newcommand{\la}{\longrightarrow}
\newcommand{\bfl}{\begin{flushleft}}
\newcommand{\efl}{\end{flushleft}}

\newcommand{\eps}{\epsilon}

\newcommand{\calj}{\mathcal{J}}

\newcommand{\xr}{\xrightarrow}

 \newcommand{\bcm}{\bar \cm}

 \newcommand{\ut}{\underbar{2}}

\begin{document}  

  \title{ Floer homotopy theory, realizing chain complexes by module spectra, and  manifolds with corners}
  \author{Ralph L. Cohen \thanks{The author was partially supported by a   grant  from the NSF} \\ Department of Mathematics  \\Stanford University \\ Stanford, CA 94305 }
 \date{\today}
\maketitle  
 \begin{abstract}  In this paper we describe and continue the study begun in \cite{cjs}   
 of the   homotopy theory that underlies Floer theory. 
  In that paper  the authors addressed the question of realizing a Floer  complex as the celluar chain complex of a $CW$-spectrum or pro-spectrum, where the   attaching maps are determined by the  compactified  moduli spaces of  connecting orbits.      The basic obstructions to the existence of this realization are the smoothness of these moduli spaces,    and    the existence of compatible collections of framings of their  stable tangent bundles.     In this note we describe a generalization of this, to show that when these moduli spaces are smooth, and 
  are oriented with respect to a generalized cohomology theory $E^*$, then a Floer $E_*$-homology theory can be defined.   In doing this we describe a functorial viewpoint on how chain complexes can be realized by $E$-module spectra, generalizing the stable homotopy realization criteria given in \cite{cjs}.  Since these moduli spaces, if smooth, will be manifolds with corners,  we  give a discussion about the appropriate notion of orientations of manifolds with corners.       \end{abstract}

 \tableofcontents

 \section*{Introduction}
 In \cite{cjs}, the authors began the study of the homotopy theoretic aspects of Floer theory.  Floer theory comes in many flavors, and it was the goal of that paper to understand the common algebraic topological properties that underly them.  A Floer stable homotopy  theory of three-manifolds in the Seiberg-Witten setting was  defined and studied by Manolescu \cite{manolescu}, \cite{manolescu2}.  This in turn was  related to the stable homotopy viewpoint of Seiberg-Witten invariants of closed $4$-manifolds studied Bauer and Furuta in \cite{bauer},  \cite{bf}.  A Floer stable homotopy type in the setting of the symplectic topology of the cotangent bundle, $T^*M$,  was shown to exist by the author in \cite{cohen}, and it was calculated to be the stable homotopy type of the free loop space, $LM$.  These were the results reported on by the author at the Abel conference in Oslo.
 
 In this largely expository note, we  discuss  basic notions of Floer homotopy type, and generalize them to discuss obstructions to the existence of a Floer $E_*$- homology theory, when $E^*$ is a generalized cohomology theory.

In a ``Floer theory"  one typically  has a functional, $\ca:  \cy \to \br$, defined on an infinite dimensional manifold $\cy$,  whose critical points generate a ``Floer chain complex" $(CF_*(\ca), \p_{\ca})$.   The boundary maps in this complex are determined by the zero dimensional moduli spaces of flow lines connecting the critical points, much the same as in classical  Morse theory.   In Morse theory, however, given a Morse function $f : M \to \br$, there is an associated handlebody decomposition of the manifold. This leads to a  corresponding $CW$-complex $X_f$  with one cell for each critical point of $f$,  which is  naturally homotopy equivalent to $M$.    The  associated cellular chain complex, $(C^f_*, \p_f)$ is equal to the Morse chain complex generated by critical points,   with  boundary homomorphisms   computed by counting gradient flow lines.   The relative attaching map between cells of $X_f$, in the case when the corresponding critical points $\alpha$ and $\beta$ labeling these cells are ``successive", was shown by Franks in \cite{franks} to be given by a map of spheres, which, under the Pontrjagin-Thom construction is represented by the compact framed manifold of flow lines,  $\cm(\alpha, \beta)$,  connecting $\alpha$ to $\beta$.  In \cite{cjs} this was generalized so as not to assume that $\alpha$ and $\beta$ are successive.  Namely it was shown   that all the stable attaching maps in the complex $X_f$ are represented by the framed cobordism classes of the  compactified moduli spaces $\bcm (\gamma, \delta)$ of flow lines.  Here these moduli spaces are viewed as  compact, framed, manifolds with corners.  

A natural question, addressed in \cite{cjs}, is under what conditions, is there an underlying $CW$-complex or spectrum, or even ``prospectrum" as described there, in which the underlying cellular chain complex is exactly the Floer complex.  In particular given a Floer functional, $\ca : \cy \to \br$,  one is looking for a $CW$-(pro)-spectrum that has one cell for each critical point, and whose relative attaching maps between
critical points $\alpha$ and $\beta$  of relative index one are determined by the   zero dimensional, compact, oriented  moduli space $\cm(\alpha, \beta)$  of flows connecting orbits between them.   Moreover, one would like this ``Floer (stable) homotopy type" to have its higher attaching maps, that is to say its entire homotopy type, determined by the higher dimensional compactified moduli spaces, $\bcm (\gamma, \delta)$, in a natural way, as in the Morse theory setting.   Roughly, the obstructions to the existence of this  Floer homotopy type was shown in \cite{cjs} to be  the existence  of the structures of smooth, framed manifolds with corners, on these compact moduli spaces.  This was expressed in terms of the notion of  ``compact, smooth, framed" categories.  

The goal of the present paper is to review and generalize these ideas, in order to address the notion of when a  ``Floer $E_*$-theory" exists,  where $E^*$- is a multiplicative, generalized cohomology theory.   That is, $E^*$ is represented by a commutative ring spectrum $E$.   (By a  ``commutative" ring spectrum, we mean an $E_\infty$-ring spectrum in the setting of algebra spectra over operads \cite{ekmm},  or  a   commutative symmetric ring spectrum as in \cite{hss}.)  We  will show  that if the compact moduli spaces  $\bcm (\gamma, \delta)$ are smooth manifolds with corners, and have a compatible family of $E^*$-orientations, then a Floer $E^*$-theory exists.  

\med 
In order to state this result more precisely, we observe that this Floer theory realization question is a special case of the question of the realization of a chain complex by an $E$-module spectrum.   When $E$ is the sphere spectrum  $ S^0$, then this question  becomes that of  the realization of a chain complex by a stable homotopy type.  To be more precise,  consider a connective chain complex  of free abelian groups,

$$
\to \cdots \to C_i \xr{\p_i} C_{i-1} \xr{\p_{i-1}} \cdots \to C_0
$$
  where we are given a basis, $\cb_i$ for the chain group $C_i$.  For example, if this is a Morse complex or a Floer complex, the bases $\cb_i$ are  given  as  the critical points of index $i$.  Now consider the tensor product  
chain complex $(C_* \otimes E_*, \p \otimes 1)$, where $E_* =  E^*(S^0)$ is the coefficient ring.   Notice that for each $i$,  one can take the free $E$-module spectrum generated by $\cb_i$,   $$ 
\ce_i = \bigvee_{\alpha \in \cb_i}  E, $$   and there is a natural isomorphism of $E_*$-modules,
$$
\pi_*(  \ce_i)  \simeq C_i \otimes E_*.
$$

\begin{definition}\label{realize}  We say that an $E$-module spectrum $X$  realizes  the complex $C_* \otimes E_*$,  if there exists a filtration 
 of $E$-module spectra converging to $X$,
 $$
 X_0 \hk X_1 \hk \cdots \hk X_k \hk \cdots X
 $$
satisfying the following properties: 
\begin{enumerate} 
\item  There is an equivalence of $E$-module spectra  of the subquotients  
$$
X_i/X_{i-1} \simeq \Sigma^i \ce_i,
$$ and 
\item
the induced composition  map in homotopy groups,
$$
C_i \otimes E_* \cong \pi_{*+i}(X_i/X_{i-1}) \xr{\delta_i}  \pi_{*+i}(\Sigma X_{i-1}) \xr{\rho_{i-1}} \pi_{*+i}(\Sigma (X_{i-1}/X_{i-2})) \cong C_{i-1}\otimes E_*
$$
is the boundary homomorphism, $\p_i \otimes 1$.     \notag
\end{enumerate}
Here the subquotient $X_j/X_{j-1}$ refers to the homotopy cofiber of the map $X_{i-1} \to X_i$, the map $\rho_j : X_j \to X_j/X_{j-1}$ is the projection map, and  the map $\delta_i : X_i/X_{i-1} \to \Sigma X_{i-1}$ is the Puppe extension of the homotopy cofibration sequence $X_{i-1} \to X_i \xr{\rho_i} X_i/X_{i-1}$.  
\end{definition}

 In \cite{cjs}, the authors introduced  a category $\calj_0$, whose objects are the nonnegative  integers $\bz^+$,  and whose space of  morphisms from $i$ to $j$ for $i > j+1$  is homeomorphic to the one point compactification of the manifold with corners,
 $$
 Mor_{\cj}(i, j) \cong (\br_+)^{i-j-1}\cup \infty
 $$
 and $Mor_{\cj}(j+1, j) = S^0$.  
 Here  $\br_+$ is the space of nonnegative real numbers.  The following   is a generalization of the realization theorem proved in \cite{cjs}.

 \begin{theorem}\label{realmod}   Let  $E-mod$  be the category  of $E$-module spectra.  Then realizations of the  chain complex $C_* \otimes E_*$ by    $E$-module spectra $X$, correspond to extensions of the association  
  $j \to \ce_j$  to functors $Z_X : \cj \to E-mod$,  with the property that for each $j$, the map
\begin{align}
  (Z_X)_{j+1, j} : Mor_{\cj}(j+1, j) &\wedge \ce_{j+1} \to \ce_{j}  \notag  \\
  S^0 &\wedge \ce_{j+1} \to \ce_{j} \notag
  \end{align}
   induces the boundary homomorphism $\p_{j+1} \otimes 1$ on the level of homotopy groups.
  \end{theorem}
  
  \med
  This theorem was  proved  for $E = S^0$ in \cite{cjs}, by displaying an explicit geometric realization of such a functor.  In this note
  we give  indicate how that construction can be extended    to prove this more general theorem.

  \med
  The next question  that we address in this paper, is how to give  geometric conditions
 on the moduli spaces of flow lines, that  induces a natural   $E$-module realization of  a Floer complex.    That is, we want to describe  geometric conditions that will induce   a functor $Z : \cj \to E-mod$ as in Theorem \ref{realmod}.   Roughly,
  the condition is that the moduli spaces are smooth, and admit a compatible family of $E^*$-orientations.
 
 The compatibility conditions of these orientations are described in terms of the ``flow-category" of a Floer functional (see \cite{cjs}).   One of the conditions required is that this category be a ``smooth, compact, category".  This notion was defined in \cite{cjs}.  Such a category is a topological category $\cc$ whose objects form a discrete set, and whose whose morphism spaces, $Mor (a,b)$ are compact, smooth manifolds with corners,  such that the composition maps $\mu_ : Mor (a,b) \times Mor (b, c) \to Mor(a, c)$ are smooth codimension one embeddings whose images lie in    the boundary.  In \cite{cohen} this was elaborated upon by describing a ``Morse-Smale" condition on such a category,  where the objects are equipped with a partial ordering, and have the notion of
 ``index" assigned to them.   The ``flow category" of a Morse function $f : M \to \br$   satisfies this condition  if the metric on $M$ is chosen so that the  Morse-Smale transversality condition  holds.  This is the condition which states that each   unstable manifold and stable manifold intersect each other transversally.   When the flow category of a Floer functional satisfies these properties, we call it a ``smooth Floer theory"  (see Definition \ref{smooth} below.)
 
  In section 3 we  study the issue of orientations of the manifolds with corners comprising the morphism spaces of a smooth, compact category.    This will involve a discussion
  of some basic properties of manifolds with corners, following the exposition given in \cite{laures}.  In particular we will make use of the fact that such  manifolds with corners 
  define  a diagrams of spaces, via the corner structure.  Our notion of an $E^*$-orientation of a smooth compact category will be  given by  a functorial collection
  of $E^*$-Thom classes of  the Thom spectra of the stable normal bundles of the morphism  manifolds.  These Thom spectra themselves are corresponding diagrams of spectra, and the orientation maps representing the Thom classes are  morphisms  of diagrams of spectra.     See Definition \ref{eorient} below.

  \med
  Our main result of this section is stated in Theorem \ref{main} below.  It says that if a smooth Floer theory  has an $E^*$- orientation of its flow category, then the Floer
  complex has a natural realization by an   $E$-module spectrum.  The homotopy groups of this spectrum can then be viewed as the ``Floer  $E_*$-homology".

 \section{Floer homotopy theory}
 
 \subsection{Preliminaries from Morse theory}
 
 In classical Morse theory  one begins with a smooth, closed $n$-manifold $M^n$, and a   smooth function $f : M^n \to \br$  with only nondegenerate critical points. Given a Riemannian metric on $M$,  one studies the flow of the gradient vector field $\nabla f$.
 In particular a flow line is a curve $\gamma : \br \to M$ satisfying the ordinary differential
 equation,
 $$
 \frac{d}{dt}\gamma(s) + \nabla f (\gamma (s)) = 0.
 $$
 By the existence and uniqueness theorem for solutions to ODE's, one knows that if
 $x \in M$ is any point, there is a unique flow line $\gamma_x$ satisfying $\gamma_x(0) = x$.  One then studies unstable and stable manifolds of the critical points,
 \begin{align}
 W^u(a) &= \{ x \in M \, : \, \lim_{t\to -\infty}\gamma_x (t) = a \}  \notag \\
 W^s(a) &= \{ x \in M \, : \, \lim_{t\to +\infty}\gamma_x (t) = a \}.   \notag
 \end{align}
 
The unstable manifold $W^u(a)$ is diffeomorphic to a disk $D^{\mu (a)}$, where $\mu (a)$ is the index of the critical point $a$.  Similarly    the stable manifold $W^s(a)$  is  diffeomorphic to a disk $D^{n-\mu (a)}$.  

For a generic choice of Riemannian metric, the unstable manifolds and stable manifolds
intersect transversally, and their intersections,
$$
W(a,b) = W^u(a) \cap W^s (b)
$$
are smooth manifolds   of  dimension equal to the relative index,  $\mu (a) - \mu (b)$.  When the choice of metric satisfies these transversality properties,  the metric is  called ``Morse-Smale".     The manifolds $W(a,b)$ have   free $\br$-actions  defined by  ``going with the flow".  That is,   for $t \in \br$, and $x \in M$, 
$$
t\cdot x = \gamma_x (t).
$$
The ``moduli space of flow lines" is the manifold
$$
\cm (a,b) = W(a,b) /\br
$$
and has dimension $\mu (a) - \mu (b) -1$.  These moduli spaces are not generally compact, but they have canonical compactifications in the following  way.

 In the case of  a Morse -Smale metric, (which we assume throughout the rest of this section),     there is a partial order on the finite set of critical points,  where $a \geq b$ if $\cm(a,b) \neq \emptyset$.   We then have
\begin{equation} \label{compact}
\bar \cm (a,b) = \bigcup_{a=a_1 >a_2> \cdots > a_k = b} \cm(a_1, a_2) \times \cdots \times \cm (a_{k-1}, a_k),
\end{equation}

The topology of $\bar \cm (a,b)$ can be described as follows.  
Since the Morse function $f : M \to \br$ is strictly decreasing along flow lines,  a flow $\gamma: \br \to M$  with the property that $\gamma (0)    \in  W(a,b)$  determines a diffeomorphism $\br \cong (f(b), f(a))$ given by the composition,
$$
 \br  \xr{\gamma}M \xr{f}  \br.
 $$
 This defines a parameterization of any $\gamma \in \cm(a,b)$ as a map
 $$
\gamma : [ f(b), f(a)] \to M
$$
that satisfies the differential equation 
\begin{equation}\label{diffeq}
\frac{d\gamma}{ds} =  \frac{\nabla f (\gamma (s))}{|f(\gamma (s))|^2},
\end{equation}
as well as the boundary conditions \begin{equation}\label{bound}
\gamma (f(b)) = b  \quad \text{and} \quad \gamma (f(a))=a.
\end{equation}  From this viewpoint, the \sl compactification \rm $\bcm(a,b)$ can be described as the space of all continuous maps $ [f(b), f(a)] \to M$ that are piecewise
smooth, (and indeed smooth off of the critical values of $f$ that lie between $f(b)$ and $f(a)$), that satisfy the differential equation (\ref{diffeq}) subject to the boundary conditions (\ref{bound}).  It is topologized with the compact open topology. 

These moduli spaces have natural framings on their stable normal bundles (or equivalently, their stable tangent bundles) in the following manner.  Let $a > b$ be critical points.  let $\eps > 0$ be chosen so that there are no critical values in the half open interval $[f(a)-\eps, f(a))$.  Define the \sl unstable sphere \rm to be the level set of the stable manifold,
$$
 S^u(a) = W^u(a) \cap f^{-1}(f(a)-\eps).
 $$
$ S^u(a)$ is a sphere of dimension $ \mu(a) - 1$. Notice there is a natural diffeomorphism,
$$
\cm(a,b) \cong S^u(a) \cap W^s(b).
$$
This leads to the following diagram,
\begin{equation}\label{intersect}
\begin{CD}
W^s(b)   @>\hk >>  M \\
@A \cup AA   @AA\cup A \\
\cm (a,b)  @>>\hk >  S^u(a).
\end{CD}
\end{equation}
From this diagram one sees that the normal bundle $\nu$ of the embedding $\cm(a,b) \hk S^u(a)$ is the restriction of the normal bundle of $W^s(b) \hk M$.  Since $W^s(b)$ is a disk, and therefore contractible, this bundle is trivial. Indeed an orientation of $W^s(b)$,  determines a homotopy class of trivialization, or a framing.  In fact  this framing determines a diffeomorphism of the bundle  to the product, $W^s(b) \times W^u(b)$.  Thus these orientations give the moduli spaces 
$\cm(a,b)$   canonical normal framings, $\nu \cong \cm(a,b) \times W^u(b) . $   

In the case when $a$ and $b$ are \sl successive \rm critical points, that is when there are no intermediate critical points $c$ with $a > c > b$, then (\ref{compact}) tells us that $\cm(a,b)$ is already compact.  This means that $\cm(a,b)$ is a closed, framed manifold, and its framed cobordism class is represented by the composition,
$$
\tau_{a,b}: S^u(a) \xr{\tau} \cm^\nu (a,b) \simeq \cm(a,b)_+ \wedge (W^u(b)\cup \infty) \xr{proj} W^u(b)\cup \infty
$$
Here $\tau : S^u(a) \to \cm^\nu(a,b)$ is the ``Thom collapse map" from the ambient sphere to the Thom space of the normal bundle.     This map between spheres was shown by Franks in \cite{franks} to be the relative attaching map in the $CW$-structure of the complex $X_f$.  

Notice that in the special case  when $\mu (a) = \mu (b) + 1$  the spheres $S^u(a)$ and $W^u(b)\cup \infty$ have the same dimension ($= \mu (b)$), and so the homotopy class of $\tau_{a,b}$ is given by its degree, which by standard considerations is the ``count" $\# \cm (a,b)$ of the number of points in the compact, oriented (framed) zero-manifold $\cm (a,b)$ (counted with sign).  Thus the ``Morse-Smale chain complex", $(C^f_*, \p)$ which is the cellular chain complex of the $CW$-complex $X_f$,  is given as follows: 

\begin{equation}\label{morsecomp}
 \to \cdots \to C^f_i \xr{\p_i} C^f_{i-1} \xr{\p_{i-1}} \cdots \to C^f_0
 \end{equation}
where $C^f_i$ is the free abelian group generated by the critical points of $f$ of index $i$, and 
\begin{equation}
\p_i [a] = \sum_{\mu (b) = i-1} \#\cm (a,b) \,  [b].
\end{equation}

The framings on the higher dimensional moduli spaces, $\cm (a,b)$, extend to their compactifications, $\bcm (a,b)$, and they have the structure of compact, framed manifolds with corners.  In \cite{cjs} it was shown that these framed manifolds define
the stable attaching maps of the the $CW$-complex $X_f$, thus generalizing Franks' theorem.   

\med
The moduli spaces $\bcm (a,b)$ fit together to define the \sl flow category,  \rm $\cc_f$,
of the Morse function $f : M \to \br$.  The objects of this  category  are the critical points of $f$, and the morphisms are given by the spaces $\bcm (a,b)$.  The compositions are the inclusions into the boundary,
$$
\bcm (a,b) \times \bcm(b,c) \hk \bcm(a,c).
$$
What the above observations describe, is the structure on the flow category $\cc_f$ of a 
\sl smooth, compact,  Morse-Smale category \rm as defined in \cite{cjs} and \cite{cohen}:

\begin{definition}\label{cptcat} A smooth, compact  category  is a topological category $\cc$ whose objects form a discrete set, and whose whose morphism spaces, $Mor (a,b)$ are compact, smooth manifolds with corners,  such that the composition maps $\mu_ : Mor (a,b) \times Mor (b, c) \to Mor(a, c)$ are smooth codimension one embeddings (of manifolds with corners) whose images lie in  the boundary.  

A smooth, compact category $\cc$  is said to be a ``Morse-Smale" category if the following
additional properties are satisfied.  
\begin{enumerate}
\item The objects of $\cc$ are partially ordered by condition
$$
 a \geq  b   \quad \text{if} \quad Mor (a,b) \neq \emptyset.
$$
\item $Mor(a,a) = \{identity \}$. 
\item  There is a set map, $\mu : Ob (\cc) \to \bz$ which preserves the partial ordering
so that  if $a > b$, 
$$
dim \, Mor(a,b) = \mu (a) - \mu (b) -1.
$$
$\mu$ is known as an ``index" map.   A Morse-Smale category such as this  is said to have finite type, if for each pair of objects  $a > b$,  there are only finitely many objects
$\alpha$ with $a > \alpha > b$. 
\end{enumerate}
\end{definition}

As was described in \cite{cjs} and \cite{cohen}, the Morse framings of the moduli spaces
$\cm (a,b)$ are  compatible with the composition structure in the flow category $\cc_f$, which gives the category the structure of a ``smooth, compact, framed" category.  We will recall this notion in section 3, when we generalize it to the notion of a ``smooth, compact, $E^*$-oriented" category, where $E^*$ is a multiplicative generalized cohomology theory.   
 It was also shown in \cite{cjs}, that the framing on this category was precisely what was needed to (functorially) realize the Morse chain complex (\ref{morsecomp}) by a stable homotopy type, which in this case was shown to be the suspension spectrum of the manifold, $\Sigma^\infty (M_+)$.  It was shown also, that if the framing was changed, one produces a change in realization, typically to the Thom spectrum of a virtual bundle, $M^\zeta$.

 \med
 \subsection{Smooth Floer theories}
In Floer theory, many, although not all of the above constructions go through.  
 In such a theory, one typically starts with a smooth map
 $$
 \ca : \cl \to \br
 $$
 where often $\cl$ is an infinite dimensional manifold, and  the functional $\ca$ is defined
 from geometric considerations. For example in symplectic Floer theory, $\cl$ is the loop space of a symplectic manifold, and $\ca$ is     the ``symplectic action" of a loop.  In ``instanton"-Floer theory,  $\cl$ is the space of gauge equivalence classes of $SU(n)$-connections on the trivial bundle over a fixed $3$-manifold $Y$, and $\ca$ is the
 Chern-Simons functional. 
 
 Often times the functional $\ca$ is then perturbed so that the critical points are isolated and nondegenerate.  Also,  a choice of metric on $\cl$ is typically shown to exist on $\ca$  so that the Morse-Smale transversality conditions hold.   So even though it is  often the case that the unstable and stable manifolds, $W^u(a)$ and $W^s(a)$ are infinite dimensional, in such a Floer theory  the intersections 
 $$
 W(a,b) = W^u(a) \cap W^s(b)
 $$
 are finite dimensional smooth manifolds.     Thus even when one cannot make sense of
 the index, $\mu (a)$, one \sl can \rm make sense of the relative index, $$\mu(a,b) = dim \, W(a,b). $$
 
 One can then form the moduli spaces, $\cm (a,b) = W(a,b)/\br$ with dimension $\mu (a,b) - 1$.  One can also form the space of ``piecewise flows" , $\bcm(a,b)$ using the same method as above (\ref{compact}). However in general these space \sl may not be compact. \rm  This is typically due to  bubbling phenomena.  Moreover, the smoothness of  $\bcm(a,b)$ may be difficult to establish.    In general this is a  deep analytic question, which has been addressed in particular examples throughout the literature
 (e.g  \cite{floer}, \cite{taubes}, \cite{donaldson},  \cite{baraudcornea}).  Furthermore a   general analytical theory of ``polyfolds"  is being developed by Hofer \cite{hofer}  to deal with these types of questions in general. However in Floer theory one generally knows that when the relative index $\mu (a,b) = 1$, then the spaces $\cm(a,b)$ is a compact, oriented, zero dimensional manifolds.  By picking a ``base" critical point $a_0$,   one can  then  can form a ``Floer chain complex" (relative to $a_0$),

 \begin{equation}\label{floercomp}
  \to \cdots \to C^{\ca}_i \xr{\p_i} C^{\ca}_{i-1} \xr{\p_{i-1}} \cdots \to C^{\ca}_0
 \end{equation}
where $C^{\ca}_i$ is the free abelian group generated by the critical points $a$ of $\ca$ of index $\mu (a, a_0) = i$, and 
\begin{equation}
\p_i [a] = \sum_{\mu (b) = i-1} \#\cm (a,b) \,  [b].
\end{equation}
By allowing the choice of base critical point $a_0$ to vary, one obtains an inverse system of Floer chain complexes as in \cite{cjs}.  Another way of handling the lack of well defined
index is to use coefficients in a  Novikov ring as in \cite{mcduffsal}.

The specific question   addressed in \cite{cjs} is the realizability of the Floer complex
by a stable homotopy type: a spectrum in the case of a fixed connective chain complex as (\ref{floercomp}),  or a prospectrum in the case of an inverse system of such chain complexes.  In this
note we  ask the following generalization of this question.  Given a commutative ring spectrum $E$,   when does  the Floer complex $C^{\ca}_* \otimes E_*$     have a  realization  by an $E$-module spectrum as in Definition (\ref{realize}) in the introduction?

Ultimately the criteria for such realizations rests on the properties of the flow category of the Floer functional, $\cc_{\ca}$.  
This is defined as in the Morse theory case,  where the objects are critical points, and the morphisms are the spaces of piecewise flow lines,  $\bcm(a,b)$.     The first criteria are that these moduli spaces can be given be given the structure of smooth, compact, manifolds with corners.  As mentioned above, this involves the analytic issues of transversality, compactness, and gluing:

\med
\begin{definition}\label{smooth}We say that a Floer functional $\ca : \cl \to \br$ generates a ``smooth Floer   theory"  if the flow category $\cc_{\ca}$ has the structure of a smooth, compact, Morse-Smale category of finite type,  as in Definition \ref{cptcat}.
\end{definition}  

\med
\noindent
\bf Remark. \rm Notice that in this definition we are assuming that a global notion
of index can be defined, as opposed to only the relative index with respect to a base critical point.   If this is not the case, one can  often replace the flow category with an inverse system of categories, each of which has a global notion of index. See \cite{cjs} for a more thorough  discussion of this point.    

\med
 Examples of  ``smooth Floer theories"  were given in \cite{cjs} and \cite{cohen}.  In particular the symplectic Floer theory of the cotangent bundle $T^*M$ was shown to have this property in \cite{cohen}.
 
In the next section we establish a general functorial description of realizations of chain complexes by $E$-module spectra,  and in section 3, we interpret these criteria geometrically, in terms of the flow category $\cc_{\ca}$ of a smooth  Floer   theory.

 \section{Realizing chain complexes by $E$-module spectra}
 
 Suppose one is given a chain complex of free abelian groups, 
\begin{equation}\label{chains}
 \to \cdots C_i \xr{\p_i} C_{i-1} \to \cdots  \to C_0
\end{equation}
 where $C_i$ has a given basis, $\cb_i$.  In the examples of a smooth  Floer theory  and Morse theory, $\cb_i$ represents the set of critical points of index $i$.  In the category of spectra,  the wedge $Z_i = \bigvee_{\beta \in \cb_i}S^0$,  where $S^o$  is the sphere spectrum,  realizes the chains  in the sense that
 $$
 H_*(Z_i) \cong C_i.
 $$
Equivalently, if $H$ represents the integral Eilenberg-MacLane spectrum, the wedge $\ch_i =  \bigvee_{\beta \in \cb_i}H$ has the property that $\pi_*(\ch_i) \cong C_i$.   More generally, given a commutative ring spectrum $E$, the wedge, 
  $\ce_i = \bigvee_{\beta \in \cb_i}E$ is a free $E$-module spectrum with the property that its homotopy groups,   
  $$
  \pi_*(\ce_i) \cong C_i \otimes E_*.
  $$
 
 The key in the realization of the chain complex $(C_* \otimes E_*, \p \otimes 1)$ by an $E$-module spectrum is to understand the role of the attaching maps.   To do this, we ``work backwards", in the sense that we study the attaching maps in an $E$-module spectrum.
 For ease of exposition, we consider a finite $E$-module spectrum.  That is,  we assume $X$ is an $E$-module that has a filtration by $E$-spectra,
 $$
X_0 \hk X_1 \hk \cdots \hk X_n = X,
$$
where each $X_{i-1} \hk X_i$ is a cofibration with cofiber, $K_i = X_i \cup c(X_{i-1})$ a free $E$-module, in that there is an equivalence
$$
K_i \simeq \bigvee_{\cdd_i}\Sigma^i E 
$$
where $\cdd_i$ is a finite indexing set.

Following the ideas of \cite{cjs},  then one can ``rebuild" the homotopy type of the $n$-fold suspension, $\Sigma^n X,$
as the union of iterated cones and suspensions of the $K_i$'s,
\begin{equation}\label{decomp}
\Sigma^n X \simeq \Sigma^nK_0 \cup c(\Sigma^{n-1} K_1) \cup \cdots \cup c^i(\Sigma^{n-i} K_i) \cup \cdots \cup c^n K_n.
\end{equation}

This decomposition can be described as follows.
Define a map $\delta_i : \Sigma^{n-i}K_i \to \Sigma^{n-i+1}K_{i-1}$ to be the iterated suspension of the composition,
 $$
 \delta_i : K_i \to \Sigma K_{i-1} \to \Sigma K_{i-1}
 $$
 where the two maps in this composition come from the cofibration sequence, $X_{i-1}\to X_i \to K_i \to \Sigma X_{i-1} \cdots$.  
$\delta_i$ is a map of $E$-module spectra.  Like was pointed out in \cite{cohen}, this induces a ``homotopy chain complex",

  \begin{equation}\label{htpychain}
K_n \xr{\delta_n} \Sigma K_{n-1} \xr{\delta_{n-1}} \cdots \xr{\delta_{i+1}}\Sigma^{n-i }K_{i } \xr{\delta_{i}} \Sigma^{n-i+1}K_{i-1} \xr{\delta_{i-1}} \cdots \xr{\delta_1}\Sigma^{n}K_0 = \Sigma^n X_0.
\end{equation}
 We refer to this as a homotopy chain complex, because examination of the defining cofibrations lead to canonical null homotopies of the compositions,  
 $$
\delta_j \circ \delta_{j+1}.
$$ 
 This canonical null homotopy
defines an extension of $\delta_j$ to the    mapping cone of $\delta_{j+1}$: 
$$
c(\Sigma^{n-j-1}K_{j+1}) \cup_{\delta_{j+1}} \Sigma^{n-j}K_j   \la \Sigma^{n-j+1}K_{j-1}.
$$
More generally for every $q$, using these null homotopies,   we have an extension  to the iterated mapping cone,
\begin{equation}\label{attach}
c^q(\Sigma^{n-j-q}K_{j+q}) \cup c^{q-1}(\Sigma^{n-j-q+1}K_{j+q-1}) \cup \cdots \cup c(\Sigma^{n-j-1}K_{j+1}) \cup_{\delta_{j+1}} \Sigma^{n-j}K_j  \la  \Sigma^{n-j+1}K_{j-1}. 
\end{equation}

In other words, for each $p > q$, these null homotopies define a map of free $E$-spectra, 
\begin{equation}\label{phi}
\phi_{p,q} :  c^{p- q-1}\Sigma^{n-p}K_p \to \Sigma^{n-q}K_q.
\end{equation}

To keep track of the combinatorics of these attaching maps,  a category $\calj$ was introduced in \cite{cjs} and used again in \cite{cohen}.   The objects of $\calj$ are the integers, $\bz$.   To describe the morphisms, we first introduced spaces $J(n,m)$ for any pair of integers $n > m$.  These spaces were defined by:   
\begin{equation}
J(n,m) = \{t_i, \, i \in \bz, \quad \text{where each $t_i$ is a nonnegative real number, and}   \quad t_i=0,  
 \, \text{unless} \, m < i < n. \}
\end{equation}
Notice that $J(n,m) \cong  \br_+^{n-m-1}$, where $\br_+^q $ is the space of $q$-tuples of   nonnegative real numbers.     Furthermore there are   natural inclusions,
$$
\iota: J(n,m) \times J(m,p)  \hk J(n,p).
$$
The image of $\iota$ consists of those sequences in $J(n,p)$ which have a $0$ in the $m^{th}$ coordinate.
This then allowed the definition of the   morphisms   in $\calj$  as follows.  For integers $n<m$ there are no morphisms from $n$ to $m$.  The only morphism from an integer $n$ to itself is the identity. If $n = m+1$,   the morphisms is defined to be the two    point space, $ Mor(m+1, m) = S^0$.   If       $n>m+1$,  $Mor(n,m)$  is given by the one point compactification,
$$
Mor (n,m) =  J(n,m)^+ = J(n,m) \cup \infty.
$$
For consistency of notation we refer to all the morphism spaces  $Mor(n,m)$ as $J(n,m)^+$.  
Composition in the category is given by addition of sequences,
$$
J(n,m)^+ \times J(m,p)^+ \to J(n,p)^+.
$$ 
 A key feature of this category is that   for a based space or spectrum $Y$,  the smash product $J(n,m)^+\wedge Y$  is  homeomorphic  iterated cone,
$$
J(n,m)^+\wedge Y = c^{n-m-1}(Y).
$$

Given integers $p > q$, then there are subcategories  $\calj^p_q$ defined to be  the full subcategory generated by integers $q \geq m \geq p.$  The category $\calj_q$ is the full subcategory of $\calj$ generated by all integers $m \geq q$.  
  
  \med
        As described in \cite{cjs},  given   a functor to the category of spaces,
 $  Z : \calj_q \to  Spaces_*$  one can take its geometric realization,
\begin{equation}\label{zrealize}
 |Z| = \coprod_{q\leq j  } Z(j)\wedge J(j, q-1)^+ / \sim
\end{equation}
 where one identifies the image of $Z(j) \wedge  J(j,i)^+\wedge J(i, q-1)^+$  in $Z(j) \wedge J(j, q-1)^+$  with its image in $Z(i) \wedge J(i, q-1)$ under the map on morphisms.  
 
 For a functor whose value is in $E$-modules, we replace the above construction of the
 geometric realization $|Z|$, by a coequalizer, in the following way:
 
 Let $Z : \calj_q \to E-mod$.  Define two maps of $E$-modules,
 \begin{equation}\label{iotamu}
 \iota, \mu :  \bigvee_{q \leq j} Z(j) \wedge  J(j,i)^+\wedge J(i, q-1)^+  \la  \bigvee_{q\leq j  } Z(j)\wedge J(j, q-1)^+.
\end{equation}
 The first map $\iota$ is induced by the inclusion, $J(j,i)\times J(i, q-1) \hk J(j, q-1).$  the second map $\mu$ is the given by the wedge of maps,  
 $$
  Z(j) \wedge J(j,i)^+\wedge J(i, q-1)^+  \xr{\mu_q \wedge 1} Z(i)\wedge J(i, q-1)^+
 $$ where $\mu_q : Z(j) \wedge J(j,i)^+ \to Z(i)$ is the action of the morphisms.
 
 \med
 \begin{definition}\label{georeal}  Given a functor $Z:  \calj_q \to E-mod$ we define its geometric realization to be the  homotopy  coequalizer (in the category $E-mod$) of the two maps,
 $$
  \iota, \mu :  \bigvee_{q \leq j} Z(j) \wedge  J(j,i)^+\wedge J(i, q-1)^+  \la  \bigvee_{q\leq j  } Z(j)\wedge J(j, q-1)^+.
  $$
  \end{definition}

  \med
  The following  is a more precise version of Theorem \ref{realmod} of the introduction.

   \med
 \begin{theorem}\label{realizing} Let $E$ be a commutative ring spectrum,  and let $$
\to \cdots  C_n \xr{\p_n} C_{n-1} \to \cdots \to C_i \xr{\p_i} C_{i-1} \xr{\p_{i-1}} \cdots \xr{\p_1} C_0,
 $$  be a chain complex of free abelian groups where $C_i$ has basis $\cb_i$.  Then
 each realization of $C_* \otimes E_*$ 
 by an $E$-module spectrum $X$,   occurs as the geometric realization $|Z_X|$  of a functor  
   $$
 Z_X :  \calj_0 \to  E-mod
 $$
 with the following properties.
 \begin{enumerate}
 \item  $Z_X(i) = \bigvee_{\cb_i} E$ \\
 \item  On morphisms of the form $J(i, i-1)^+ = S^0$,  the functor $Z_X:  \bigvee_{\cb_i} E \to  \bigvee_{\cb_{i-1}} E$  induces the boundary  homomorphism on homotopy groups:  
 $$\p_i : C_i \otimes E_* \to C_{i-1} \otimes E_*.$$
 \end{enumerate}
 \end{theorem}

\med
\begin{proof}  The proof mirrors  the arguments given in \cite{cjs} and \cite{cohen} which concern realizing chain complexes of abelian groups by stable homotopy types.  We therefore present a sketch of the constructions and indicate the modifications needed to prove this theorem.

First suppose that $Z :  \calj_0 \to E-mod$ is a functor satisfying the properties (1) and (2) as  stated in the theorem.  The geometric realization
$|Z|$ has a natural filtration by $E$-module spectra,
$$
|Z|_0 \hk |Z|_1 \hk \cdots \hk |Z|_{k-1} \hk |Z|_k \hk \cdots  |Z|
$$
where the $k^{th}$-filtration, $|Z|_k$, is the homotopy coequalizer of the maps $\iota$ and $\mu$ as in Definition \ref{georeal}, but  the wedges   only involve the first $k$-terms.  Notice that the subquotients (homotopy cofibers)  are homotopy equivalent to a wedge,
$$
 |Z|_k/|Z|_{k-1}  \simeq \bigvee_{\cb_k}\Sigma^k E \simeq \Sigma^k Z(i).
 $$
The second equivalence holds by property (1) in the theorem. 
Furthermore  by property 2,     the composition defining the attaching map, 
$$
C_k \otimes E_* \cong \pi_{*+k}(|Z|_k/|Z|_{k-1} )  \to  \pi_{*+k}(\Sigma |Z|_{k-1}) \to \pi_{*+k}(\Sigma (|Z|_{k-1}/|Z|_{k-2} ) \cong C_{k-1}\otimes E_*
$$
is the boundary homomorphism, $\p_k \otimes 1$.  Thus the geometric realization of such a functor $Z : \calj_0 \to E-mod$ gives a realization of the
chain complex $C_* \otimes E_*$ according to Definition \ref{realize}.

Conversely, suppose that the filtered $E-module$
 $$
 X_0 \hk X_1 \hk \cdots \hk X_k \hk \cdots X
 $$
 realizes $C_*\otimes E_*$ as in Definition \ref{realize}.  We define a functor
 $$
 Z_X : \calj_0 \to E-mod
 $$
 by letting $Z_X (i) = \bigvee_{\cb_i} E = \Sigma^{-i}X_i/X_{i-1}.$  To define the value of $Z_X$ on morphisms, we need to define
 compatible maps of $E-modules$,
\begin{align}
 J(i,j)^+ \wedge Z_X(i) &\to Z_X(j)  \notag \\
  c^{i-j-1} \Sigma^{-i}X_i/X_{i-1}  &\to \Sigma^{-j}X_{j}/X_{j-1} \notag
  \end{align}
  
 These are  defined to be the maps $\phi_{i,j}$ given by    the attaching maps in a filtered $E$-module as  described above (\ref{phi}).
 \end{proof}

 \section{Manifolds with corners,  $E^*$-orientations of flow categories, and     Floer $E_*$ - homology }
 
 Our goal in this section is to apply Theorem \ref{realizing}  to a Floer chain complex, in order to determine when a smooth Floer theory (Definition \ref{smooth}) can be realized by
 an $E$-module, so as to produce a ``Floer $E_*$- homology theory".  To do this we want to identify the appropriate properties of the moduli spaces $\bcm(a,b)$  defining the morphisms in the flow category $\cc_{\ca}$ so that they induce a functor $Z_\ca : \calj_0 \to E-mod$ that satisfies Theorem \ref{realizing}.  The basic structure these moduli spaces need to possess is an  $E^*$- orientation on their stable normal bundles. The functor will then be defined via a Pontrjagin-Thom construction.  
 
  Now these moduli spaces have a natural  corner structure, and these $E^*$-orientations must be compatible with this corner structure.  Therefore we begin this section with a general discussion of manifolds
 with corners, their embeddings, and induced normal structures.    This general discussion follows that of \cite{laures}.  
 
  \med
 Recall that an $n$-dimensional manifold with corners $M$  has charts which are local homeomorphisms with  $\br_+^n$.    Let $\psi : U \to (\br_+)^n$ be a chart of a manifold with corners $M$.  For $x \in U$,  the number of zeros of this chart, $c(x)$ is independent of the chart.   We define a \sl face \rm of $M$ to be a connected component of the space  $\{m{\in}M \, \text{such that} \, c(m)=1\}$.
 
 Given an integer $k$, there is a notion of a manifold with corners having ``codimension k",  or a $<k>$-manifold.  We recall the definition from \cite{laures}. 
 
 \begin{definition}\label{kmanifold} A $<k>$-manifold  is a manifold  with corners   $M$  together with an ordered $k$-tuple  $(\partial{_1}M,...,\partial{_k}M)$ of unions of faces of $M$  satisfying the following properties.
 \begin{enumerate}
 \item Each $m{\in}M$ belongs to $c(m)$ faces 
 \item $\partial{_1}M {\cup} {\cdots} {\cup} \partial_{k} M = \partial {M} $ 
 \item For all $1{\leq}i{\neq}j{\leq}k$, $\partial{_i}M {\cap} \partial{_j}M$ is a face of both  $\partial{_i}M$ and $\partial{_j}M$.
 \end{enumerate}
 
\end{definition}
 
  The  archetypical example of a  $\langle k \rangle$-manifold is $ \br_+^k$.  In this case
 the face $F_j \subset \br_+^k$ consists of those $k$-tuples with the $j^{th}$- coordinate equal to zero.

As described in \cite{laures}, the data of a $\langle k \rangle$-manifold can be encoded in a   categorical way  as follows.  Let $\underbar{2}$ be the partially ordered set with two objects, $\{0, 1\}$, generated by a single nonidentity morphism $0 \to 1$.  Let $\ut^k$ be the product of $k$-copies of the category $\ut$.  A $\langle k \rangle$-manfold $M$ then defines a functor from $\ut^k$   to the category of topological spaces,  where for an object  $a = (a_1, \cdots , a_k) \in \ut^k$,
$M(a)$ is the intersection of the faces $\p_i(M)$   with $a_i = 0$.  Such a functor is a $k$-dimensional cubical diagram of spaces, which, following Laures' terminology, we refer to as a $\langle k \rangle$-diagram, or a $\langle k \rangle-space$.    Notice that $\br_+^k(a) \subset \br_+^k$ consists of those $k$-tuples of nonnegative real numbers so that the $i^{th}$-coordinate is zero for every $i$ such that $a_i=0$.  More generally, consider the $<k>$-Euclidean space, $\br_+^k \times \br^n$,  where the value on $a \in \ut^k$ is $\br_+^k(a) \times \br^n$.      In general  we refer to a functor $\phi : \ut^k \to \cc$ as a $<k>$-object in the category $\cc$.

    In this section we will consider embeddings of manifolds with corners into Euclidean spaces $M \hk \br_+^k \times \br^n$ of the form given by the following definition.

\begin{definition}\label{embed} A ``neat embedding" of a $\langle k \rangle$-manifold $M$ into Euclidean space $ \br_+^k \times \br^m  $  is a natural transformation of $\langle k \rangle$-diagrams
$$
e : M \hk  \br_+^k \times \br^m
$$   that satisfies the following properties:
\begin{enumerate}
\item  For each $a{\in}\underline{2}^k$, $e(a)$ is an embedding.
\item
For all $b < a$, the intersection   $M(a) \cap  \left( \mathbb{R}^k_+(b)  \times  \mathbb{R}^m\right) = M(b)$,  and this intersection is perpendicular.  That is,    there is some $\epsilon > 0$ such that $$M(a) \cap  \left(\mathbb{R}^k_+(b)  \times  [0,\epsilon)^k(a-b) \times \mathbb{R}^m\right) = M(b)  \times [0,\epsilon)^k(a-b).$$  
\end{enumerate}
Here $a-b$ denotes the object of $\ut^k$ obtained by subtracting the $k$-vector $b$ from the $k$-vector $a$. 
\end{definition}

\med
In \cite{laures} it was proved   that every $<$k$>$-manifold neatly embeds in $\mathbb{R}^{k}_+{\times}R^N$ for $N$ sufficiently large.  In fact it was proved there that a manifold with corners $M$  admits  a neat embedding into $\mathbb{R}^{k}_+{\times}R^N$ \sl if and only if \rm $M$ has the structure of a $<k>$-manifold.  Furthermore in \cite{genauer} it was proved that the connectivity of the space of neat embeddings, $Emb_{<k>}(M; \mathbb{R}^{k}_+{\times}R^N)$ increases with the dimension $N$.

\med
Notice that given  an embedding of manifolds with corners,  $e : M \hk  \br_+^k \times \br^m$,  it has a well defined normal bundle.  In particular, for any pair of objects in $\ut^k$ $a > b$, then the normal bundle of $e(a) : M(a) \hk \br^m \times \br_+^k(a)$, when restricted to $M(b)$, is the normal bundle of $e(b) : M(b) \hk \br^m \times \br_+^k(b)$. 

\med

Now  let  $\cc$  be the flow category of a smooth Floer theory.  Recall from Definition \ref{smooth}  that this means  $\cc$  is a smooth,  compact, Morse-Smale category.     By assumption, the morphism spaces, $Mor(a,b) = \bcm (a,b)$    are compact, smooth, manifolds with corners.  Notice that in fact they are
$\langle  k(a,b)\rangle$ - manifolds, where $k(a,b) = \mu (a) - \mu (b) -1$.  This structure is given by  the faces
 $$
\p_j (\bcm (a,b)) = \bigcup_{\mu (c) = \mu(a) - j} \bcm(a, c) \times \bcm(c,b).
$$
 These faces clearly satisfy the  intersection property necessary for being a $\langle k(a,b) \rangle$- manifold (Definition \ref{kmanifold}).   The condition on the flow category $\cc$ necessary for it to induce a geometric realization of the Floer complex by $E$-modules, will be that the moduli spaces, $\bcm (a,b)$ admit
 a compatible family of $E^*$-orientations.  To define this notion carefully, we first observe that  a  commutative ring spectrum $E$   induces a $<k>$-diagram in the category of spectra (``$<k>$ - spectrum"),  $E<k>$, defined in the following recursive manner.
 
 For $k = 1$, we let $E<1>: \ut \to Spectra$  be defined by $E<k>(0) = S^0$, the sphere spectrum, and $E<1>(1) = E$.  The image of the morphism $0 \to 1$ is the unit of the ring spectrum $S^0 \to E$.  
 
 Now recursively, suppose we have defined the $<j>$-spectrum $E<j>$ for every $1\leq j \leq k-1$.  We now define $E<k>$.    Let $\iota_k = (1,1, \cdots ,1)$ be the maximal element
 in the poset $\ut^k$.  We define $E<k>(\iota_k) = E$.  Now consider the faces, 
 $  \p_i (\iota_k) = (1, 1, \cdots , 0, 1, \cdots, 1)$ with the single $0$ in this vector occurring in the $i^{th}$ slot. We can think of this element as $(\iota_{i-1}, 0, \iota_{k-i})$,
 when viewed as the image of the natural inclusion functor $$\ut^{i-1} \times \{0\} \times \ut^{k-i} \hk \ut^k.$$
The image  consists of all sequences that have a $0$ in the $i^{th}$ slot.

  We then define
\begin{align}
 E<k>(\p_i \iota_k) &= E<i-1>(\iota_{i-1}) \wedge S^0 \wedge E<k-i>(\iota_{k-i}) \notag \\
 &= E \wedge S^0 \wedge E. \notag
 \end{align}
 The image of the morphism $\p_i \iota_k \hk \iota_k$  is given by the ring multiplication,
 $$
  E \wedge S^0 \wedge E \to E.
  $$
 Notice that every other object of $a \in \ut^k$ is less than or equal to $\p_i \iota_k$ for some $i = 1, \cdots , k$, and so we can identify $a$ as an element of the product category  $\ut^{i-1}\times {0} \times \ut^{k-i} \hk \ut^k$.   In other words, $a = (a_{i-1}, 0, a_{k-i})$ for some
 $a_{i-1} \in \ut^{i-1}$,  $a_{k-i}\in \ut^{k-i}.$   We then define
 $$
 E<k>(a) = E<i-1>(a_{k-i}) \wedge S^0 \wedge E<k-i>(a_{k-i}).
 $$
 By the associativity of the ring multiplication, the above data suffices to give a well defined functor,
 $$
 E<k> : \ut^k \to Spectra.
 $$
 Furthermore the ring structure induces pairings
 $$
m_{k,r}: E<k> \wedge E<r> \to E<k+r>
 $$
 which are natural transformations of functors $\ut^{k+r} \to Spectra$  (i.e maps
 of $<k+r>$-spectra).
 
 \med
 This structure allows us to define one more construction.  Suppose $\cc$ is a smooth,  compact,  Morse-Smale category of finite type as in Definition \ref{cptcat}.  We can then define an associated category,  $E_{\cc}$ whose objects are the same as the objects of $\cc$, and whose morphisms are given by the spectra,
 $$
 Mor_{E_{\cc}}(a,b)  = E<k(a,b)>
 $$
 where $k(a,b) =  \mu(a) - \mu (b) -1$.  Here  $\mu (a)$ is the index of the object $a$ as in Definition \ref{cptcat}.  The composition law is the pairing,
\begin{align}
 E<k(a,b)> \wedge E<k(b,c)> &=  E<k(a,b)> \wedge S^0 \wedge  E<k(b,c)>  \notag\\
 &\xr{1 \wedge u \wedge 1} E<k(a,b)> \wedge E<1> \wedge  E<k(b,c)>  \notag \\
 &\xr{\mu} E<k(a,c)>. \notag
 \end{align}
  Here $u : S^0 \to E = E<1>$ is the unit.    This category encodes the multiplication in the ring spectrum $E$.  
  
  \med
  We now make precise what  it means to say that the moduli spaces in a Floer flow category, $\bcm (a,b)$, or more generally, the morphisms in a smooth, compact, Morse-Smale category,   have an $E^*$-orientation. 
  
 \med
 Let $M$ be a $<k>$-manifold, and let $e : M \hk \br_+^k \times \br^N$ be a neat embedding.  The Thom space, $Th (M, e)$ has the structure of an $<k>$-space, where for $a \in \ut^k$, $Th (M, e)(a)$ is the Thom space of the normal bundle of the associated embedding, $M(a) \hk \br_+^k (a) \times \br^N$.  We can then desuspend, and define the Thom spectrum, $M^\nu_e = \Sigma^{-N}Th (M,e)$ to be the associated $<k>$-spectrum. The Pontrjagin-Thom construction defines a map of $<k>-$ spaces, $$\tau_e : \left( \br_+^k \times \br^N\right)\cup \infty = (( \br_+^k)\cup \infty) \wedge S^N \to  Th (M, e).$$ Desuspending we get a map of $<k>$-spectra, $\Sigma^\infty (( \br_+^k)\cup \infty) \to  M^\nu_e. $  Notice that  the homotopy type (as  $<k>$-spectra) of $M^\nu_e$ is independent of the embedding $e$. We denote the homotopy type of this normal Thom spectrum as $M^\nu$, and the homotopy type of the Pontrjagin-Thom map,  $\tau : \Sigma^\infty (( \br_+^k)\cup \infty) \to M^\nu$.
 
 We define an $E^*$-normal orientation to be a cohomology class (Thom class), represented by a map $u : M^\nu \to E$ such that cup product defines an isomorphism,
 $$
 \cup u : E^*(M) \xr{\cong} E^*(M^\nu).
 $$
 
 In order to solve the realization problem, we need to have coherent $E^*$-orientations of all the moduli spaces, $\bcm (a,b)$- making up a smooth, Floer theory.
 In order to make this precise, we make the following definition.
 
 \med
 \begin{definition}\label{normalthom} Let $\cc$ be the flow category of a smooth Floer theory.  Then a ``normal Thom spectrum" of the category $\cc$, is a category $\cc^\nu$, enriched over spectra, that satisfies the following properties.
 \begin{enumerate}
 \item The objects of $\cc^\nu$ are the same as the objects of $\cc$.
 \item The morphism spectra  $Mor_{\cc^\nu}(a,b)$ are $<k(a,b)>$-spectra,  having  the homotopy type of the normal Thom spectra $\bcm(a,b)^\nu$, as $<k(a,b)>$-spectra. The composition maps, $$\circ : Mor_{\cc^\nu}(a,b)\wedge Mor_{\cc^\nu}(b,c) \to Mor_{\cc^\nu}(a,c)$$
 have the homotopy type of the maps, $$
 \bcm(a,b)^\nu \wedge \bcm(b,c)^\nu \to \bcm(a,c)^\nu$$ of normal bundles corresponding to the composition maps in $\cc$, $\bcm(a,b) \times \bcm(b,c) \to \bcm (a,c)$.  Recall that these   are maps of $<k(a,c)>$-spaces induced by the inclusion of a component of the boundary.   
 \item The morphism spectra are equipped with   maps  $\tau_{a,b} : \Sigma^\infty (J(\mu (a), \mu (b))^+) = \Sigma^\infty ( (  \br_+^{k(a,b)})\cup \infty)) \to Mor_{\cc^\nu}(a,b)$  that are of the homotopy type of    the  Pontrjagin-Thom  $\tau:  \Sigma^\infty ( (  \br_+^{k(a,b)})\cup \infty)) \to \bcm(a,b)^\nu$, and so that the following diagrams commute:
 $$
 \begin{CD}
   \Sigma^\infty (J(\mu (a), \mu (b))^+) \wedge  \Sigma^\infty (J(\mu (b), \mu (c))^+)  @>>>  \Sigma^\infty (J(\mu (a), \mu (c))^+)  \\
   @V\tau_{a.b}\wedge \tau_{b,c} VV    @VV\tau_{a,c} V \\
   Mor_{\cc^\nu}(a,b) \wedge Mor_{\cc^\nu}(b,c)  @>>> Mor_{\cc^\nu}(a,c).
   \end{CD}
   $$ Here the top horizontal map is defined via the composition maps  in the category $\calj$,  and the bottom horizontal map is defined via the composition maps in $\cc^\nu$.  
   \end{enumerate}
   \end{definition}
   
   Notice that according to this definition, a normal Thom spectrum  is simply a functorial collection of Thom spectra of stable normal bundles.  Finite, smooth, compact categories have such normal Thom spectra, defined via embeddings into Euclidean spaces (see \cite{cohen}).  
   
   \med
   We can now define what it means to have coherent $E^*$-orientations on the moduli spaces $\bcm (a,b)$.  
   
   \med
   \begin{definition}\label{eorient}
   An $E^*$-orientation of a flow category of a smooth Floer theory, $\cc$, is a functor, $u : \cc^\nu \to  E_{\cc}$, where $\cc^\nu$ is a normal Thom spectrum of $\cc$, such that on morphism spaces, the  induced map
   $$
   Mor_{\cc^\nu}(a,b) \to   E<k(a,b)>
   $$
   is a map of $<k(a,b)>$-spectra that defines an $E^*$ orientation of $ Mor_{\cc^\nu}(a,b) \simeq \bcm(a,b)^\nu$. 
   \end{definition}
   
   Notice that when $E = S^0$, the sphere spectrum, then an $E^*$-orientation as defined here, is equivalent to a \sl framing \rm of the category $\cc$, as defined in \cite{cjs}, and  \cite{cohen}.  
   
   \med
   We can now indicate the proof of the main theorem of this section.
   
   \med
   \begin{theorem}\label{main}.  Let $\cc$ be the flow category of a smooth Floer theory defined by a functional $\ca : \cl \to \br$, and let $E$ be a commutative ring spectrum. Let $u : \cc^\nu \to E_{\cc}$ be an $E^*$-orientation of the category $\cc$.   Then the induced Floer complex $(\cc^{\ca}_*, \p)$ has a natural realization by an $E$-module $|Z^{\ca}_E|$.
   \end{theorem}
   
   \noindent
   \bf Remarks. \rm  
   \begin{enumerate}
   \item In the case when $E = S^0$, this is a formulation of the result in \cite{cjs} that says that a Floer theory that defines a \sl framed, \rm smooth, compact flow category defines a ``Floer stable homotopy type".  
   \item When $E = H\bz$, the integral Eilenberg-MacLane spectrum, then the notion of an $H\bz$-orientation of a flow category $\cc$ is equivalent to the existence of ordinary ``coherent
   orientations" of the moduli spaces, $\bcm (a,b)$ as was studied in \cite{floerhofer}.
   \item  Given a flow category $\cc$ satisfying the hypotheses of the theorem with respect to a ring spectrum $E$,  then one defines the homotopy groups of the 
   corresponding realization, $\pi_*(|Z^{\ca}_E|)$ to be the  ``Floer $E_*$-homology of functional $\ca$".
   \end{enumerate}

   \med
   \begin{proof}   Let $\cc$ be such a category, and $u : \cc^\nu \to E_{\cc}$ an $E^*$-orientation.  Consider the corresponding Floer complex (see (\ref{floercomp})):
   $$
   \to \cdots \to C^{\ca}_i \xr{\p_i} C^{\ca}_{i-1} \xr{\p_{i-1}} \cdots \to C^{\ca}_0
$$  Recall that the chain group  $C^{\ca}_i$ has a basis, $\cb_i$ consisting of the
objects (critical points) of index $i$, and the  boundary maps are defined by
counting (with sign) the elements of the oriented moduli spaces $\cm(a,b)$, where 
$a$ and $b$ have relative index one ($\mu (a) = \mu (b) +1.$)   

We define a functor $Z^{\ca}_E : \calj_0 \to E-mod$ that satisfies the conditions of Theorem \ref{realizing}.  As required, on objects we define
$$
Z^\ca_E (i) = \bigvee_{\cb_i}E.
$$
   
   To define $Z^\ca_E$ on the morphisms in $\calj_0$, we need to define, for every $i >j$
   $E$-module maps
\begin{align}
   Z^\ca_E (i,j):  J(i,j)^+ \wedge Z^\ca_E (i)  &\to Z^\ca_E (j) \notag \\
\bigvee_{\cb_i}  J(i,j)^+ \wedge E  &\to \bigvee_{\cb_j}E \notag
\end{align}  
That is, for each $a \in \cb_i$, and $b \in \cb_j$, we need to define an $E$-module map
$$
 Z^\ca_E (a,b) : J(i,j)^+ \wedge E \to E,
 $$
 or, equivalently, an ordinary map of spectra,
 $$
  Z^\ca_E (a,b) : \Sigma^\infty (J(i,j)^+) \to E.
  $$  This map is defined by the $E^*$-orientation of the moduli space $\bcm(a,b)$, as the composition,
 $$
   \Sigma^\infty (J(i, j)^+)  \xr{\tau_{a,b}} Mor_{\cc^\nu}(a,b)   \xr{u} E<i-j-1>
   $$
   where $\tau_{a,b}$ is the Pontrjagin-Thom collapse map of the normal Thom spectrum
  $Mor_{\cc^\nu}(a,b) \simeq \bcm(a,b)^\nu$ (see Defintion \ref{normalthom}), and $u$ is the map   of $<i-j-1>$-spectra coming from the $E^*$-orientation  (Definition \ref{eorient}). The functorial properties of these maps is precisely what is required
  so that they fit together to define a functor,
 $Z^{\ca}_E : \calj_0 \to E-mod$.   It is straightforward to check that this  functor satisfies 
  the conditions of Theorem \ref{realizing}, and so by that theorem, $|Z^{\ca}_E|$ is an $E$-module spectrum that realizes the Floer complex  $C^{\ca}_* \otimes E_*$.
    \end{proof}


\begin{thebibliography}{99}{ 




\bibitem{baraudcornea} J.-F Barraud and O. Cornea,  \emph{Lagrangian intersections and the Serre spectral sequence},  preprint (2004),  arXiv: mathDG/0401094

\bibitem{bauer} S. Bauer, \emph{A stable cohomotopy refinement of Seiberg-Witten invariants:II}, Invent. Math. \textbf{155(1)}, (2004) 21-40.

\bibitem{bf}S. Bauer and M. Furuta, \emph{A stable cohomotopy refinement of Seiberg-Witten invariants: I},  Invent. Math. \textbf{155(1)}, (2004) 1-19






\bibitem{cohen}  R.L. Cohen, \emph{The Floer homotopy type of the cotangent bundle}, preprint (2007),  arXiv:math/0702852.





\bibitem{cjs}R.L. Cohen, J.D.S. Jones, and G.B. Segal, \emph{Floer's infinite dimensional Morse theory and homotopy theory}
  \sl Floer Memorial Volume, \rm
Birkhauser Verlag Prog.
in Math. \bf vol. 133 \rm  (1995), 287 - 325.

 

\bibitem{cjscat} R.L. Cohen, J.D.S. Jones, and G.B. Segal, \emph{Morse theory and classifying spaces}, Stanford University preprint, (1995)  available at http://math.stanford.edu/\~{}ralph/papers.html 

\bibitem{donaldson}S. K. Donaldson,   \bf Floer Homology groups in Yang-Mills theory \rm  with the assistance of M. Furuta and D. Kotschick,  \sl Cambridge Tracts in Mathematics \bf 147, \rm Cambridge Univ. Press, (2002).  


\bibitem{ekmm} A.~D.~Elmendorf, I.~Kriz, M.~A.~Mandell, J.~P.~May, {\em Rings, modules, and algebras in stable homotopy theory. {W}ith an appendix by M.~Cole}, Mathematical Surveys and Monographs, {\bf 47}, American Mathematical Society, Providence, RI, 1997. 


\bibitem{floer} A. Floer,  \emph{Morse theory for Lagrangian intersections} J. Differential Geometry \bf 28 \rm (1988),  513-547.


\bibitem{floerhofer} A. Floer and H. Hofer, \emph{Coherent orientations for periodic orbit problems in symplectic geometry}, Math. Zeit., \bf 212 \rm (1993), 13-38. 

\bibitem{franks} J.M.  Franks,  \emph{Morse-Smale flows and homotopy theory}, Topology \bf{18} \rm no. 3, (1979), 199-215.



\bibitem{genauer} J. Genauer,  \emph{Stanford University PhD thesis}, in preparation.


\bibitem{hofer} H. Hofer, \emph{A generaized Fredholm theory and applications}   	In: Current Developments in Mathematics,
D. Jerison et al. (eds.), International Press (2006)  math.SG/0509366

\bibitem{hss} M.~Hovey, B.~Shipley, J.~Smith, {\em Symmetric spectra}, J.\ Amer.\ Math.\ Soc.\ {\bf 13} (2000), 149--208. 


 
 
 
\bibitem{laures}G. Laures, \emph{On cobordism of manifolds with corners}, Trans. of AMS \bf 352 \rm no. 12,  (2000), 5667-5688.

\bibitem{manolescu} C. Manolescu, \emph{Seiberg-Witten-Floer stable homotopy type of three-manifolds with $b_1=0$}, Geom. Topol. \textbf{7}, (2003), 889-932

\bibitem{manolescu2}C. Manolescu, \emph{ A gluing theorem for the relative Bauer-Furuta invariants} , J. Diff. Geom. Topol. \textbf{76}, (2007), 889-932

\bibitem{mcduffsal} D. McDuff and D. Salamon, \textbf{$J$-holomorphic curves and Quantum cohomology}, American Mathematical Society,  Providence, (1994).




\bibitem{salamon} D. Salamon, \emph{Lectures on Floer homology}, Symplectic geometry and topology (Y. Eliashberg and L. Traynor, eds.), IAS/Park City Mathematics Series, AMS (1999), 143-225.


      
     
     \bibitem{taubes}  C. Taubes, $L^2$\bf - moduli spaces on $4$-manifolds with cylindrical  ends,  \rm \sl Monographs in Geometry and Topology, I, \rm International Press (1993).

    

      }\end{thebibliography}
 \end{document}